\documentclass[a4paper,11pt]{amsart}
\usepackage{amsfonts,stmaryrd,amssymb,amsmath,a4wide,verbatim,url,color,epsfig,mathrsfs,graphicx,pgf,amsthm,amsfonts,mathtext,cite,enumerate,float,mathtools,tikz}
\usepackage[british]{babel}
\usepackage{hyperref}
\usetikzlibrary{arrows}
\usepackage[matrix,arrow,curve]{xy}
\usepackage[outline]{contour}
\usepackage{mathtools}
\usepackage{tcolorbox}
\usepackage{letltxmacro}
\usepackage{nameref}
\usepackage{lipsum}
\usepackage{enumitem}
\usepackage{xcolor}
\usepackage{subfig}
\makeatletter
\@addtoreset{equation}{section}\makeatother

\theoremstyle{plain}

\theoremstyle{definition}

\theoremstyle{definition}

\makeatletter
\newcommand*{\rom}[1]{\expandafter\@slowromancap\romannumeral #1@}
\makeatother

\begin{document}

\captionsetup[subfigure]{labelformat=empty}

\title[Cosmology of Plane Geometry]{Cosmology of Plane Geometry}



\author[Alexander Skutin]{Alexander Skutin}

\maketitle

\section{Introduction}

This paper focuses on a new approach to plane geometry and develops important concepts that can allow researchers to unite and observe plane geometry from a new, meaningful perspective. The present short note is a first chapter of the paper, referenced in the last section.


\section{Deformation principle}

In the area of plane geometry, the deformation principle refers to replacing a certain configuration of points, lines, and circles with more general ones, in which equal points, lines, or circles in the original (undeformed) configuration are replaced by their deformed versions, i.e. points, lines, and circles that are not equal in the general case but are related to one another. If one is aware that points, lines, or circles are equal in the case of non-deformed configuration, one can predict the deformed versions’ connections in terms of general configuration. Consideration of the non-deformed case aids in predicting which points should be connected in the general deformed case.

\subsection{Basic deformation example}

For the first example of application of deformation principle consider the configuration of a square $ABCD$ with its center $O$. Now, we can look on the point $O$ as on the third vertices $O_{ab}$, $O_{bd}$, $O_{cd}$, $O_{da}$ of the triangles $O_{ab}AB$, $O_{bc}BC$, $O_{cd}CD$, $O_{da}DA$ constructed internally on the sides of $ABCD$ which are isosceles right angled triangles ($\angle AO_{ab}B = 90^{\circ}$, $|O_{ab}A| = |O_{ab}B|$ and so on). So, in the case when $ABCD$ is a square we see that $O=O_{ab}=O_{bc}=O_{cd}=O_{da}$ are equal, thus, the deformation principle predicts that in the general (deformed) case of an arbitrary $ABCD$, deformed points $O_{ab}$, $O_{bc}$, $O_{cd}$, $O_{da}$ should be connected. And in fact, they are: $O_{ab}O_{cd}$ is perpendicular to $O_{bc}O_{da}$ and they have equal lengths.

\medskip

Therefore, we can state the next result

\medskip

$\textbf{Theorem 1 (Deformation of a square with its center).}$ Consider any (convex)\\ quadrilateral $ABCD$ and let points $O_{ab}$, $O_{bc}$, $O_{cd}$, $O_{da}$ are chosen inside of $ABCD$ such that $O_{ab}AB$, $O_{bc}BC$, $O_{cd}CD$, $O_{da}DA$ are isosceles right angled triangles (see picture below).

\medskip

Then the following properties of the points $O_{ab}$, $O_{bc}$, $O_{cd}$, $O_{da}$ are true:

\medskip

\begin{enumerate}
    \item $O_{ab}O_{cd}\perp O_{bc}O_{da}$
    \item $|O_{ab}O_{cd}| = |O_{bc}O_{da}|$.
\end{enumerate}


In the previous theorem we see how the deformation principle works if we deform a square with its center, similarly we can deform other simple configurations and get complicated results. Each time when we deform we should have decided which points we like to deform and how we deform them. For example, in the previous case we can look on the square's center $O$ as on the intersection of the internal angle $\angle ABC$, $\angle BCD$, $\angle CDA$, $\angle DAB$ bisectors which will provide us with another type of deformation of $ABCDO$ (in the case of the angle bisectors the resulting deformational fact may be stated as : Points $O_1$, $O_2$, $O_3$, $O_4$ are cyclic). To see why the deformation principle is natural we can look on it as on the reverse process to the process of consideration of particular cases of theorems. For example, Theorem 1 becomes trivial in the particular case of a square because in this case points $O_{ab}$, $O_{bc}$, $O_{cd}$, $O_{da}$ become equal. So, the particular case should be seen as an evidence of the general one.

\subsection{Possible deformations}

Theorem 1 is the deformation of four nested points into four connected points in the general case. Similarly, other objects equal (or trivially connected) in the non-deformed case can be deformed to understand their connected deformations. The following table shows which objects may be considered non-deformed and how their deformations may appear:

\bigskip

\begin{table}[htb]
\begin{tabular}{l|l}
$\textbf{Undeformed objects}$ & $\textbf{Possible deformations}$                                                                                                                                                                                         \\ \hline
Coincide points    & Collinear points, concyclic points                                                                                                                                                                                   \\ \hline
Coincide lines     & Concurrent lines                                                                                                                                                                                                 \\ \hline
Coincide circles   & \begin{tabular}[c]{@{}l@{}}Concurrent circles, coaxial circles, circles having\\ radical line which has many nice properties wrt\\ original configuration (if there are two circles)\end{tabular} \\ \hline
Coincide triangles & \begin{tabular}[c]{@{}l@{}}Perspective triangles, triangles which vertices\\ are lying on the same conic\end{tabular}                                                                                                \\ \hline
Concyclic points   & \begin{tabular}[c]{@{}l@{}}Concyclic points, points lying on the same conic\\ \end{tabular}    \\ \hline                                                                                                                      
\end{tabular}


\end{table}

\hspace{-6mm}

\subsection{Deformations of an equilateral triangle}

Instead of a square $ABCD$, this paper uses an equilateral triangle $ABC$ and attempts to deform its center, incircle, circumcircle, and many other closely related objects. Deformation of an equilateral triangle is the most powerful tool for producing and verifying results throughout this paper.

As an example of the deformation of an equilateral triangle, we can consider an equilateral triangle with center $O$ and interpret $O$ as the entire set of Clark Kimberling centers $X_i$ (see \cite{Cl}). Therefore, the deformation theory predicts that, in the general case of an arbitrary triangle $ABC$, the triangle centers $X_i$ should have many relations among themselves, such as collinearity, concyclity, etc. Many of these relationships may be found in \cite{Cl}.

Next we introduce some first examples of application of the deformation principle in the case of equilateral triangles.


\textbf{Example 1.} Consider any triangle $ABC$. Let $A'BC$ be the equilateral triangle constructed on the side $BC$ such that $A'$, $A$ are on the same half plane wrt $BC$. Denote the center $O_a$ of $A'BC$. Similarly define $O_b$, $O_c$. In the case of equilateral $ABC$ we get that $O_a=O_b=O_c$~-- coincides with the center $O$ of $ABC$. So, we can predict that in the general case of an arbitrary $ABC$ we will have that the points $O_a$, $O_b$, $O_c$ will be connected and have relations with the base triangle $ABC$. And, in fact, the next such relation can be formulated: Points $O_a$, $O_b$, $O_c$ form the equilateral triangle which circumcircle is passing through the first Fermat point of $ABC$.

\textbf{Example 2.} Consider a triangle $ABC$ with the first and second Fermat points $F_1$, $F_2$. Let $F_a$, $F_b$, $F_c$ be the second Fermat points of $F_1BC$, $F_1AC$, $F_1AB$, respectively. In the case of equilateral $ABC$ we get that $F_1$ coincides with the center $O$ of $ABC$ and $F_2$ coincides with a point $P$ lying on the circumcircle of $ABC$. Aditionally, in the equilateral triangle case we have that $F_a$, $F_b$, $F_c$ are the reflections of center $O$ of $ABC$ wrt its sides and the points $A$, $B$, $C$, $F_a$, $F_b$, $F_c$ form a regular hexagon which circumcircle contains $P$. Thus, we can predict that in the general case of an arbitrary $ABC$ the configuration $ABCF_aF_bF_cF_2$ will inherit some properties of a regular hexagon and a point $P$ on its circumcircle. And, in fact, the next such relation can be formulated: Point $F_2$ lies on the circumcircle of $F_aF_bF_c$.

\textbf{Example 3.} Consider a triangle $ABC$ and any point $P$. Let $AP$, $BP$, $CP$ meet the circumcircle $(ABC)$ of $ABC$ second time at $A'$, $B'$, $C'$ (i.e. $A'B'C'$ is the circumcevian triangle of $P$ wrt $ABC$). Consider the nine-point centers $N$, $N_a$, $N_b$, $N_c$ of $ABC$, $A'BC$, $B'AC$, $C'AB$, respectively. Denote

\begin{enumerate}
    \item the reflections $N_a'$, $N_b'$, $N_c'$ of $N_a$, $N_b$, $N_c$ wrt $BC$, $AC$, $AB$, respectively
    \item the reflections $N_a''$, $N_b''$, $N_c''$ of $N_a$, $N_b$, $N_c$ wrt midpoints of $BC$, $AC$, $AB$, respectively.
\end{enumerate}

In the case of equilateral $ABC$ and $P=O$~-- its center we get that $N_a'=N_b'=N_c'=N=P=O$ and $N_a''=N_b''=N_c''=N=P=O$. So, we can predict that in the general case of an arbitrary $ABC$ and an point $P$ we will have that the points $N_a'$, $N_b'$, $N_c'$, $N$, $P$ ($N_a''$, $N_b''$, $N_c''$, $N$, $P$) will be connected with each other. And, in fact, the next such connections can be formulated:

\begin{enumerate}
    \item Points $N_a'$, $N_b'$, $N_c'$, $N$ lie on the same circle
    \item Points $N_a''$, $N_b''$, $N_c''$, $N$ lie on the same circle.
    
\end{enumerate}

\section{Complete article}

The full version of this article was written by the Author involving geometric problems contribution by Tran Quang Hung, Kadir Altintas and Antreas Hatzipolakis.

Download it from the links below

     Cosmology of Plane Geometry: Concepts and Theorems (2019) \href{https://www.scribd.com/document/421475794}{Scribd/421475794}
    
     Cosmology of Plane Geometry (Improved version 2021) \href{https://www.scribd.com/document/510674976}{Scribd/510674976}.

\end{document}